\documentclass[reqno,12pt]{amsart}
\usepackage{amssymb}
\usepackage{amsmath}
\usepackage{amsfonts}

\newtheorem{theorem}{Th\'eor\`eme}[section]
\theoremstyle{plain}

\newtheorem{corollary}[theorem]{Corollaire}

\newtheorem{definition}[theorem]{D\'efinition}

\newtheorem{lemma}[theorem]{Lemme}

\newtheorem{remark}[theorem]{Remarque}

\numberwithin{equation}{section}
\input{tcilatex}

\begin{document}
\title{ Caract\'{e}risation d'une norme Fr\'{e}chet diff\'{e}rentiable pour
les Banach duaux complexes }
\author{Daher Mohammad}
\address{117, Rue Gustave Courbet 77350 Le M\'{e}e Sur Seine-France}
\email{daher.mohammad@ymail.com}

\begin{abstract}
Let $X$ be a complex Banach space; in this work we characterize the property
of Fr\'{e}chet diffrentiability for the dual space $X^{\ast }.$
\end{abstract}

\maketitle

0\bigskip AMS Classification : 46B70

Mots cl\'{e}s : Interpolation, lisse\ \ \ \ \ \ \ \ \ \ \ \ \ \ \ \ \ \ \ \
\ \ \ \ \ \ \ \ \ \ \ \ \ \ \ \ \ \ \ \ \ \ \ \ \ \ \ \ \ \ \ \ \ \ \ \ \ \
\ \ \ \ \ \ \ \ \ \ \ \ \ \ \ \ \ \ \ \ \ \ \ \ \ \ \ \ \ \ \ \ \ \ \ \ \ \
\ \ \ \ \ \ \ \ \ \ \ \ \ \ \ \ \ \ \ \ \ \ \ \ \ \ \ \ \ \ \ \ \ \ \ \ \ \
\ \ \ \ \ \ \ \ \ \ \ \ \ \ \ \ \ \ \ \ \ \ \ \ \ \ \ \ \ \ \ \ \ \ \ \ \ \
\ \ \ \ \ \ \ \ \ \ \ \ \ \ \ \ \ \ \ \ \ \ \ \ \ \ \ \ \ \ \ \ \ \ \ \ \ \
\ \ \ \ \ \ \ \ \ \ \ \ \ \ \ \ \ \ \ \ \ \ \ \ \ \ \ \ \ \ \ \ \ \ \ \ \ \
\ \ \ \ \ \ \ \ \ \ \ \ \ \ \ \ \ \ \ \ \ \ \ \ \ \ \ \ \ \ \ \ \ \ \ \ \ \
\ \ \ \ \ \ \ \ \ \ \ \ \ \ \ \ \ \ \ \ \ \ \ \ \ \ \ \ \ \ \ \ \ \ \ \ \ \
\ \ \ \ \ \ \ \ \ \ \ \ \ \ \ \ \ \ \ \ \ \ \ \ \ \ \ \ \ \ \ \ \ \ \ \ \ \
\ \ \ \ \ \ \ \ \ \ \ \ \ \ \ \ \ \ \ \ \ \ \ \ \ \ \ \ \ \ \ \ \ \ \ \ \ \
\ \ \ \ \ \ \ \ \ \ \ \ \ \ \ \ \ \ \ \ \ \ \ \ \ \ \ \ \ \ \ \ \ \ \ \ \ \
\ \ \ \ \ \ \ \ \ \ \ \ \ \ \ \ \ \ \ \ \ \ \ \ \ \ \ \ \ \ \ \ \ \ \ \ \ \
\ \ \ \ \ \ \ \ \ \ \ \ \ \ \ \ \ \ \ \ \ \ \ \ \ \ \ \ \ \ \ \ \ \ \ \ \ \
\ \ \ \ \ \ \ \ \ \ \ \ \ \ \ \ \ \ \ \ \ \ \ \ \ \ \ \ \ \ \ \ \ \ \ \ \ \
\ \ \ \ \ \ \ \ \ \ \ \ \ \ \ \ \ \ \ \ \ \ \ \ \ \ \ \ \ \ \ \ \ \ \ \ \ \
\ \ \ \ \ \ \ \ \ \ \ \ \ \ \ \ \ \ \ \ \ \ \ \ \ \ \ \ \ \ \ \ \ \ \ \ \ \
\ \ \ \ \ \ \ \ \ \ \ \ \ \ \ \ \ \ \ \ \ \ \ \ \ \ \ \ \ \ \ \ \ \ \ \ \ \
\ \ \ \ \ \ \ \ \ \ \ \ \ \ \ \ \ \ \ \ \ \ \ \ \ \ \ \ \ \ \ \ \ \ \ \ \ \
\ \ \ \ \ \ \ \ \ \ \ \ \ \ \ \ \ \ \ \ \ \ \ \ \ \ \ \ \ \ \ \ \ \ \ \ \ \
\ \ \ \ \ \ \ \ \ \ \ \ \ \ \ \ \ \ \ \ \ \ \ \ \ \ \ \ \ \ \ \ \ \ \ \ \ \
\ \ \ \ \ \ \ \ \ \ \ \ \ \ \ \ \ \ \ \ \ \ \ \ \ \ \ \ \ \ \ \ \ \ \ \ \ \
\ \ \ \ \ \ \ \ \ \ \ \ \ \ \ \ \ \ \ \ \ \ \ \ \ \ \ \ \ \ \ \ \ \ \ \ \ \
\ \ \ \ \ \ \ \ \ \ \ \ \ \ \ \ \ \ \ \ \ \ \ \ \ \ \ \ \ \ \ \ \ \ \ \ \ \
\ \ \ \ \ \ \ \ \ \ \ \ \ \ \ \ \ \ \ \ \ \ \ \ \ \ \ \ \ \ \ \ \ \ \ \ \ \
\ \ \ \ \ \ \ \ \ \ \ \ \ \ \ \ \ \ \ \ \ \ \ \ \ \ \ \ \ \ \ \ \ \ \ \ \ \
\ \ \ \ \ \ \ \ \ \ \ \ \ \ \ \ \ \ \ \ \ \ \ \ \ \ \ \ \ \ \ \ \ \ \ \ \ \
\ \ \ \ \ \ \ \ \ \ \ \ \ \ \ \ \ \ \ \ \ \ \ \ \ \ \ \ \ \ \ \ \ \ \ \ \ \
\ \ \ \ \ \ \ \ \ \ \ \ \ \ \ \ \ \ \ \ \ \ \ \ \ \ \ \ \ \ \ \ \ \ \ \ \ \
\ \ \ \ \ \ \ \ \ \ \ \ \ \ \ \ \ \ \ \ \ \ \ \ \ \ \ \ \ \ \ \ \ \ \ \ \ \
\ \ \ \ \ \ \ \ \ \ \ \ \ \ \ \ \ \ \ \ \ \ \ \ \ \ \ \ \ \ \ \ \ \ \ \ \ \
\ \ \ \ \ \ \ \ \ \ \ \ \ \ \ \ \ \ \ \ \ \ \ \ \ \ \ \ \ \ \ \ \ \ \ \ \ \
\ \ \ \ \ \ \ \ \ \ \ \ \ \ \ \ \ \ \ \ \ \ \ \ \ \ \ \ \ \ \ \ \ \ \ \ \ \
\ \ \ \ \ \ \ \ \ \ \ \ \ \ \ \ \ \ \ \ \ \ \ \ \ \ \ \ \ \ \ \ \ \ \ \ \ \
\ \ \ \ \ \ \ \ \ \ \ \ \ \ \ \ \ \ \ \ \ \ \ \ \ \ \ \ \ \ \ \ \ \ \ \ \ \
\ \ \ \ \ \ \ \ \ \ \ \ \ \ \ \ \ \ \ \ \ \ \ \ \ \ \ \ \ \ \ \ \ \ \ \ \ \
\ \ \ \ \ \ \ \ \ \ \ \ \ \ \ \ \ \ \ \ \ \ \ \ \ \ \ \ \ \ \ \ \ \ \ \ \ \
\ \ \ \ \ \ \ \ \ \ \ \ \ \ \ \ \ \ \ \ \ \ \ \ \ \ \ \ \ \ \ \ \ \ \ \ \ \
\ \ \ \ \ \ \ \ \ \ \ \ \ \ \ \ \ \ \ \ \ \ \ \ \ \ \ \ \ \ \ \ \ \ \ \ \ \
\ \ \ \ \ \ \ \ \ \ \ \ \ \ \ \ \ \ \ \ \ \ \ \ \ \ \ \ \ \ \ \ \ \ \ \ \ \
\ \ \ \ \ \ \ \ \ \ \ \ \ \ \ \ \ \ \ \ \ \ \ \ \ \ \ \ \ \ \ \ \ \ \ \ \ \
\ \ \ \ \ \ \ \ \ \ \ \ \ \ \ \ \ \ \ \ \ \ \ \ \ \ \ \ \ \ \ \ \ \ 

\section{\textsc{Introduction}}

\bigskip Soit $X$ un espace de Banach; notons $B_{X}$ sa boule unit\'{e} ferm%
\'{e} et S$_{X}$ sa sph\`{e}re unit\'{e}. Notons d'autre part, $\left\langle
x,x^{\ast }\right\rangle $ l'accoumplement entre un \'{e}l\'{e}ment de $X$
et un \'{e}l\'{e}ment de $X^{\ast }$ le dual de $X.$

\begin{definition}
\label{tz}Un espace de Banach $X$ est Fr\'{e}chet diff\'{e}rentiable si,
pour tout $x\in X-\left\{ 0\right\} ,$
\end{definition}

lim$_{\left\Vert h\right\Vert \rightarrow 0}\frac{\bigl\Vert x+h\bigr\Vert+%
\bigl\Vert x-h\bigr\Vert-2\bigl\Vert x\bigr\Vert}{\bigl\Vert h\bigr\Vert}=0.$

\begin{definition}
\bigskip \label{UL}Soit $X$ un espace de Banach; $X$ est uniform\'{e}ment
lisse, si pour tout $x\in S_{X},$ il existe un unique point $x^{\ast }\in
S_{X^{\ast }}$ tel que pour tout $\varepsilon >0,$ il existe $\delta >0$ v%
\'{e}rifiant $\left\Vert x^{\ast }-g\right\Vert <\varepsilon $ d\`{e}s que $%
g\in B_{X^{\ast }}$, et $\left\langle x,g\right\rangle >1-\delta .$
\end{definition}

\bigskip Soit $X$ un espace de Banach r\'{e}el; d'apr\`{e}s le r\'{e}sultat
de \v{S}mulyan \cite[Chap.I,Th.1.4,Cor.1.5]{D-G-Z}, l'espace dual $X^{\ast }$
est Fr\'{e}chet diff\'{e}rentiable si et seulement s'il v\'{e}rifie l'une
des propri\'{e}t\'{e}s suivantes:

1-$X^{\ast }$ est uniform\'{e}ment lisse.

2- pour toutes suites $(x_{n})_{n\geq 0},(y_{n})_{n\geq 0}$ dans $S_{X}$ et
tout $a^{\ast }\in S_{X^{\ast }}$ tels que $\left\langle x_{n},a^{\ast
}\right\rangle \rightarrow _{n\rightarrow +\infty }1$ et $\left\langle
y_{n},a^{\ast }\right\rangle \rightarrow _{n\rightarrow +\infty }1,$ alors $%
\bigl\Vert x_{n}-y_{n}\bigr\Vert\rightarrow _{n\rightarrow +\infty }0.$

Dans \cite[Lemma 5.2]{Da} on montre que si le Banach X est complexe, alors $%
X^{\ast }$ est Fr\'{e}chet diff\'{e}rentiable si et seulment si $X^{\ast }$ v%
\'{e}rifie la propri\'{e}t\'{e} 2. La d\'{e}monstration dans le cas complexe
est exacetement similaire \`{a} la d\'{e}monstration de \cite[Chap.I,Th.1.4]%
{D-G-Z} avec une petite modification.

Soit $X$ un espace de Banach complexe; dans ce travail on montre que $%
X^{\ast }$ est uniform\'{e}ment lisse si et seulement si $X^{\ast }$ v\'{e}%
rifie la propri\'{e}t\'{e} 2. \ Dans le suite, on montre que si $X^{\ast }$
est G\^{a}teaux diff\'{e}rentiable, alors $L^{p}(\Omega ,X)^{\ast }$ est G%
\^{a}teaux diff\'{e}rentiable ($(\Omega ,\mu )$ est un espace de probabilit%
\'{e} et $1<p<+\infty $).

\begin{theorem}
\label{K}Soit $X$ un espace de Banach. Alors $X^{\ast }$ est uniform\'{e}%
ment lisse si et seulement si, $X^{\ast }$ v\'{e}rifie la propri\'{e}t\'{e}
2.
\end{theorem}

D\'{e}monstration: Supposons que $X^{\ast }$ est uniform\'{e}ment lisse.
Soient $(f_{n})_{n\geq 0},(g_{n})_{n\geq 0}$ deux suites dans $S_{X}$ et $%
x^{\ast }\in S_{X^{\ast }}$ tels que $\left\langle f_{n},x^{\ast
}\right\rangle \rightarrow _{n\rightarrow +\infty }1$ et $\left\langle
g_{n},x^{\ast }\right\rangle \rightarrow _{n\rightarrow +\infty }1.$ Soit
d'autre part $\varepsilon >0;$ comme $X^{\ast }$ est uniform\'{e}ment lisse,
il existe $x^{\ast \ast }$ un point unique dans $S_{X^{\ast \ast }}$ et $%
\delta >0$ v\'{e}rifiant $\bigl\Vert x^{\ast \ast }-y^{\ast \ast }\bigr\Vert<%
\frac{\varepsilon }{2}$ d\`{e}s que $y^{\ast \ast }\in B_{X^{\ast \ast }}$
et $\left\langle x^{\ast },y^{\ast \ast }\right\rangle >1-\delta .$ Notons $%
u_{n}=\frac{f_{n}}{\left\langle f_{n},x^{\ast }\right\rangle }\left\vert
\left\langle f_{n},x^{\ast }\right\rangle \right\vert $ et $v_{n}=\frac{g_{n}%
}{\left\langle g_{n},x^{\ast }\right\rangle }\left\vert \left\langle
g_{n},x^{\ast }\right\rangle \right\vert .$ Comme $\lambda _{n}=\frac{%
\left\vert \left\langle f_{n},x^{\ast }\right\rangle \right\vert }{%
\left\langle f_{n},x^{\ast }\right\rangle }\rightarrow _{n\rightarrow
+\infty }1$ et $\mu _{n}=\frac{\left\vert \left\langle g_{n},x^{\ast
}\right\rangle \right\vert }{\left\langle g_{n},x^{\ast }\right\rangle }%
\rightarrow _{n\rightarrow +\infty }1,$ $\left\langle u_{n},x^{\ast
}\right\rangle \rightarrow _{n\rightarrow +\infty }1$ et $\left\langle
v_{n},x^{\ast }\right\rangle \rightarrow _{n\rightarrow +\infty }1.$ Mais, $%
\left\langle u_{n},x^{\ast }\right\rangle ,$ $\left\langle v_{n},x^{\ast
}\right\rangle \in \mathbb{R}^{+},$ il existe donc $n_{0}$ tel que $%
\left\langle u_{n},x^{\ast }\right\rangle >1-\delta $ et $\left\langle
v_{n},x^{\ast }\right\rangle >1-\delta ,$ pour tout $n\geq n_{0};$ par cons%
\'{e}quent $\bigl\Vert x^{\ast \ast }-u_{n}\bigr\Vert<\frac{\varepsilon }{2}$
et $\bigl\Vert x^{\ast \ast }-v_{n}\bigr\Vert<\frac{\varepsilon }{2},$ ceci
implique que $\bigl\Vert u_{n}-v_{n}\bigr\Vert<\varepsilon .$ Il en r\'{e}%
sulte que $\bigl\Vert u_{n}-v_{n}\bigr\Vert\rightarrow _{n\rightarrow
+\infty }0.$ D'autre part, $\ f_{n}=\frac{u_{n}}{\lambda _{n}}$, $g_{n}=%
\frac{v_{n}}{\mu _{n}},$ $\lambda _{n}\rightarrow _{n\rightarrow +\infty }1$
et $\mu _{n}\rightarrow _{n\rightarrow +\infty }1,$ donc $\bigl\Vert %
f_{n}-g_{n}\bigr\Vert\rightarrow _{n\rightarrow +\infty }0.\blacksquare $

Supposons que $X^{\ast }$ n'est pas uniform\'{e}ment lisse. Il existe donc $%
a^{\ast }\in X^{\ast },$ $a^{\ast \ast }\in X^{\ast \ast },$ $\varepsilon >0$
et une suite $(x_{n}^{\ast \ast })_{n\geq 0}$ dans la boule unit\'{e} de $%
X^{\ast \ast }$ tels que $\left\langle a^{\ast },a^{\ast \ast }\right\rangle
=1,$ $\bigl\Vert x_{n}^{\ast \ast }-a^{\ast \ast }\bigr\Vert>\varepsilon $
et 
\begin{equation}
\left\langle a^{\ast },x_{n}^{\ast \ast }\right\rangle >1-\frac{1}{n+1}
\label{DF}
\end{equation}

Pour tout $n\geq 0,$ il existe $u_{n}^{\ast }$ dans la boule unit\'{e} de $%
X^{\ast }$ tel que%
\begin{equation}
\left\langle u_{n},x_{n}^{\ast \ast }-a^{\ast \ast }\right\rangle
>\varepsilon .  \label{LL}
\end{equation}%
Pour tout $n\geq 0$ notons%
\begin{eqnarray*}
D_{n} &=&\left\{ x^{\ast \ast }\in B_{X^{\ast \ast }};\text{ }\left\vert
\left\langle a^{\ast },x^{\ast \ast }\right\rangle -1\right\vert <\frac{1}{%
n+1}\right\} \cap \\
&&\left\{ x^{\ast \ast }\in B_{X^{\ast \ast }};\text{ }\left\vert
\left\langle u_{n}^{\ast },x^{\ast \ast }\right\rangle -\left\langle
u_{n}^{\ast },a^{\ast \ast }\right\rangle \right\vert <\frac{1}{n+1}\right\}
\end{eqnarray*}%
et%
\begin{eqnarray*}
E_{n} &=&\left\{ x^{\ast \ast }\in B_{X^{\ast \ast }};\text{ }\left\vert
\left\langle a^{\ast },x_{{}}^{\ast \ast }\right\rangle -\left\langle
a^{\ast },x_{n}^{\ast \ast }\right\rangle \right\vert <\frac{1}{n+1}\right\}
\cap \\
&&\left\{ x^{\ast \ast }\in B_{X^{\ast \ast }};\text{ }\left\vert
\left\langle u_{n}^{\ast },x_{{}}^{\ast \ast }\right\rangle -\left\langle
u_{n}^{\ast },x_{n}^{\ast \ast }\right\rangle \right\vert <\frac{1}{n+1}%
\right\} .
\end{eqnarray*}

Comme $D_{n}$ est un voisinage pr\'{e}faible de $a^{\ast \ast }$ et $E_{n}$
est un voisinage pr\'{e}faible de $x_{n}^{\ast \ast },$ il existe $x_{n}\in
D_{n}\cap X$ et $y_{n}\in E_{n}\cap X$ (car $B_{X}$ est pr\'{e}faiblement
dense dans $B_{X^{\ast \ast }})$. Donc%
\begin{equation}
\left\langle x_{n},a^{\ast }\right\rangle \rightarrow _{n\rightarrow +\infty
}1,  \label{X}
\end{equation}%
\begin{equation}
\left\langle x_{n},u_{n}^{\ast }\right\rangle -\left\langle u_{n}^{\ast
},a^{\ast \ast }\right\rangle \rightarrow _{n\rightarrow +\infty }0,
\label{B}
\end{equation}

\begin{equation}
\left\langle y_{n},a^{\ast }\right\rangle -\left\langle a^{\ast
},x_{n}^{\ast \ast }\right\rangle \rightarrow _{n\rightarrow +\infty }0
\label{KU}
\end{equation}

et%
\begin{equation}
\left\langle y_{n},u_{n}^{\ast }\right\rangle -\left\langle \left\langle
u_{n}^{\ast },x_{n}^{\ast \ast }\right\rangle \right\rangle \rightarrow
_{n\rightarrow +\infty }0.  \label{H}
\end{equation}

D'autre part, d'apr\`{e}s (\ref{DF}), $\left\langle a^{\ast },x_{n}^{\ast
\ast }\right\rangle \rightarrow _{n\rightarrow +\infty }1,$ ceci implique
d'apr\`{e}s (\ref{KU}) que%
\begin{equation}
\left\langle y_{n},a^{\ast }\right\rangle \rightarrow _{n\rightarrow +\infty
}1.  \label{O}
\end{equation}

D'apr\`{e}s (\ref{B}) et \ref{H} il existe $n_{0}$ tel que $\left\vert
\left\langle x_{n},u_{n}^{\ast }\right\rangle -\left\langle u_{n}^{\ast
},a^{\ast \ast }\right\rangle \right\vert <\frac{\varepsilon }{4}$ et $%
\left\vert \left\langle u_{n}^{\ast },x_{n}^{\ast \ast }\right\rangle
-\left\langle y_{n},u_{n}^{\ast }\right\rangle \right\vert <\frac{%
\varepsilon }{4},$ pour tout $n\geq n_{0}.$ Il en r\'{e}sulte d'apr\`{e}s (%
\ref{LL}), que 
\begin{eqnarray*}
\bigl\Vert x_{n}-y_{n}\bigr\Vert &\geq &\left\vert \left\langle
x_{n}-y_{n},u_{n}^{\ast }\right\rangle \right\vert \\
&\geq &\left\vert \left\langle u_{n}^{\ast },a^{\ast \ast }\right\rangle
-\left\langle u_{n}^{\ast },x_{n}^{\ast \ast }\right\rangle \right\vert
-\left\vert \left\langle x_{n},u_{n}^{\ast }\right\rangle -\left\langle
u_{n}^{\ast },a^{\ast \ast }\right\rangle \right\vert -\left\vert
\left\langle u_{n}^{\ast },x_{n}^{\ast \ast }\right\rangle -\left\langle
y_{n},u_{n}^{\ast }\right\rangle \right\vert \geq \frac{\varepsilon }{2},
\end{eqnarray*}

pour tout $n\geq n_{0}.$ Comme $\left\langle x_{n},a^{\ast }\right\rangle
\rightarrow _{n\rightarrow +\infty }1$ et $\left\langle y_{n},a^{\ast
}\right\rangle \rightarrow _{n\rightarrow +\infty }1$ d'apr\`{e}s (\ref{X})
et (\ref{O})$,$ $X^{\ast }$ ne v\'{e}rifie pas la propri\'{e}t\'{e} 2.$%
\blacksquare $

D'apr\`{e}s le th\'{e}or\`{e}me \ref{K} et le lemme 5.2 de \cite{Da} on a le
corollaire suivant:

\begin{corollary}
\label{wq}Soit $X$ un espace de Banach complexe. Les assertions suivantes
sont \'{e}quivalentes.
\end{corollary}

1)- $X^{\ast }$ est Fr\'{e}chet diff\'{e}rentiable.

2- pour toutes suites $(x_{n})_{n\geq 0},(y_{n})_{n\geq 0}$ dans $S_{X}$ et
tout $a^{\ast }\in S_{X^{\ast }}$ tels que $\left\langle x_{n},a^{\ast
}\right\rangle \rightarrow _{n\rightarrow +\infty }1$ et $\left\langle
y_{n},a^{\ast }\right\rangle \rightarrow _{n\rightarrow +\infty }1,$ alors $%
\bigl\Vert x_{n}-y_{n}\bigr\Vert\rightarrow _{n\rightarrow +\infty }0.$

3)- $X^{\ast }$ est uniform\'{e}ment lisse.

\begin{remark}
\label{GT}Soit $X$ un espace de Banach complexe. Alors $X$ est uniform\'{e}%
ment lisse si et seulement si, pour toutes suites $(f_{n}^{\ast })_{n\geq
0},(g_{n}^{\ast })_{n\geq 0}$ dans $S_{X^{\ast }}$ et tout $x\in S_{X}$ tels
que $\left\langle x,f_{n}^{\ast }\right\rangle \rightarrow _{n\rightarrow
+\infty }1$ et $\left\langle x,g_{n}^{\ast }\right\rangle \rightarrow
_{n\rightarrow +\infty }1,$ alors $\bigl\Vert f_{n}^{\ast }-g_{n}^{\ast }%
\bigr\Vert\rightarrow _{n\rightarrow +\infty }0.$
\end{remark}

La d\'{e}monstration de cette remarque est similaire \`{a} \cite[%
Chap.I,Th.1.4]{D-G-Z}.

\begin{definition}
\label{mmm}L'espace de Banach $(Y,\left\Vert .\right\Vert )$ est G\^{a}teaux
diff\'{e}rentiable s'il est G\^{a}teaux diff\'{e}rentiable en tout point $a$
dans $Y-\left\{ 0\right\} ,$ c' est \`{a} dire si, pour tout $h\in Y,$ lim$_{%
\mathbb{R}\ni t\rightarrow 0}\frac{\left\Vert a+th\right\Vert -\left\Vert
a\right\Vert }{t}$ existe$.$
\end{definition}

Soient $(\Omega ,\mu )$ un espace de probabilit\'{e}, $1<p<+\infty $ et $X$
un espace de Banach s\'{e}parable$;$ notons

\begin{equation*}
VB^{p}(\Omega ,X^{\ast })=\left\{ f:\Omega \rightarrow X^{\ast }\text{ pr%
\'{e}faiblement mesurable et }\mathop{\displaystyle \int}\limits_{\Omega }%
\bigl\Vert f(\omega )\bigr\Vert_{X^{\ast }}^{p}d\mu (\omega )<+\infty
\right\} .
\end{equation*}

Dans \cite{Le-Su} on montre que si le Banach X est Fr\'{e}chet diff\'{e}%
rentiable, alors $L^{p}(\Omega ,X)$ a la m\^{e}me propri\'{e}t\'{e}. Dans le
proposition suivante on montre un r\'{e}sultat analogue lorsque $X^{\ast }$
est G\^{a}teaux diff\'{e}rentiable.

\begin{lemma}
\label{vvv}Soit $X$ un espace de Banach; les assertions suivantes sont \'{e}%
quivalents :
\end{lemma}

1) $X$ est G\^{a}teaux diff\'{e}rentiable.

2) Pour tout $n\in \mathbb{N}$ et tout $x\in S_{X},$ lim$_{R\ni t\rightarrow
0}\frac{\left\Vert x+th\right\Vert ^{2^{n}}+\left\Vert x-th\right\Vert
^{2^{n}}-2\left\Vert x\right\Vert ^{2^{n}}}{t}=0,$ pour tout $h\in X.$

3) Il existe $n\in \mathbb{N}$ tel que pour tout $x\in S_{X},$ lim$_{R\ni
t\rightarrow 0}\frac{\left\Vert x+th\right\Vert ^{2^{n}}+\left\Vert
x-th\right\Vert ^{2^{n}}-2\left\Vert x\right\Vert ^{2^{n}}}{t}=0,$ pour tout 
$h\in X.$

D\'{e}monstration: 1)$\implies 2).$

\emph{Etape 1.} Soient $n\in \mathbb{N}$, $x\in S_{X}$ et $h\in X;$ montrons
que lim$_{R\ni t\rightarrow 0}\frac{\left\Vert x+th\right\Vert
^{2^{n}}-\left\Vert x\right\Vert ^{2^{n}}}{t}$ existe$.$

En effet,

Nous allons montrer l'\'{e}tape 1 par r\'{e}currence sur $n.$ L'\'{e}nonc%
\'{e} est \'{e}vidament vrai pour $n=0.$ Supposons que lim$_{R\ni
t\rightarrow 0}\frac{\left\Vert x+th\right\Vert ^{2^{n}}-\left\Vert
x\right\Vert ^{2^{n}}}{t}$ existe; montrons que lim$_{R\ni t\rightarrow 0}$ $%
\frac{\left\Vert x+th\right\Vert ^{2^{n+1}}-\left\Vert x\right\Vert
^{2^{n+1}}}{t}$ existe.

Comme $(\left\Vert x+th\right\Vert ^{2^{n+1}}-\left\Vert x\right\Vert
^{2^{n+1}})=(\left\Vert x+th\right\Vert ^{2^{n}}-\left\Vert x\right\Vert
^{2^{n}})(\left\Vert x+th\right\Vert ^{2^{n}}+\left\Vert x\right\Vert
^{2^{n}})$, lim$_{R\ni t\rightarrow 0}\frac{\left\Vert x+th\right\Vert
^{2^{n+1}}-\left\Vert x\right\Vert ^{2^{n+1}}}{t}=$lim$_{R\ni t\rightarrow 0}%
\frac{2\left\Vert x\right\Vert ^{2^{n}}\left[ \left\Vert x+th\right\Vert
^{2^{n}}-\left\Vert x\right\Vert ^{2^{n}}\right] }{t}$ existe$.$

\emph{Etape 2.} Soient $n\in \mathbb{N}$ et $h\in X.$ Montrons que lim$%
_{R\ni t\rightarrow 0}\frac{\left\Vert x+th\right\Vert ^{2^{n}}+\left\Vert
x-th\right\Vert ^{2^{n}}-2\left\Vert x\right\Vert ^{2^{n}}}{t}=0.$

En effet, d'apr\`{e}s l'\'{e}tape 1, nous avons que

\begin{eqnarray*}
&&\text{lim}_{R\ni t\rightarrow 0}\frac{\left\Vert x+th\right\Vert
^{2^{n}}+\left\Vert x-th\right\Vert ^{2^{n}}-2\left\Vert x\right\Vert
^{2^{n}}}{t} \\
&=&\text{lim}_{R\ni t\rightarrow 0}\frac{\left\Vert x+th\right\Vert
^{2^{n}}-\left\Vert x\right\Vert ^{2^{n}}}{t}+\text{lim}_{R\ni t\rightarrow
0}\frac{\left\Vert x-th\right\Vert ^{2^{n}}-\left\Vert x\right\Vert ^{2^{n}}%
}{t}=0.\blacksquare
\end{eqnarray*}

2)$\implies 3).$ \'{e}vident.

3)$\implies 1).$

Supposons qu'il existe $n\geq 1$ tel que pour tout $x\in S_{X}$ et tout $%
h\in X,$ lim$_{R\ni t\rightarrow 0}\frac{\left\Vert x+th\right\Vert
^{2^{n}}+\left\Vert x-th\right\Vert ^{2^{n}}-2\left\Vert x\right\Vert
^{2^{n}}}{t}=0.$ Il suffit de montrer que 
\begin{equation*}
\text{lim}_{R\ni t\rightarrow 0}\frac{\left\Vert x+th\right\Vert
^{2^{n-1}}+\left\Vert x-th\right\Vert ^{2^{n-1}}-2\left\Vert x\right\Vert
^{2^{n-1}}}{t}=0.
\end{equation*}

Remarquons que

\begin{eqnarray}
&&\left\Vert x+th\right\Vert ^{2^{n}}+\left\Vert x-th\right\Vert
^{2^{n}}-2\left\Vert x\right\Vert ^{2^{n}}  \notag \\
= &&\left[ (\left\Vert x+th\right\Vert ^{2^{n}}+\left\Vert x-th\right\Vert
^{2^{n}})^{\frac{1}{2}}-\sqrt{2}\left\Vert x\right\Vert ^{2^{n-1}}\right]
\label{gt} \\
&&\times \left[ (\left\Vert x+th\right\Vert ^{2^{n}}+\left\Vert
x-th\right\Vert ^{2^{n}})^{\frac{1}{2}}+\sqrt{2}\left\Vert x\right\Vert
^{2^{n-1}}\right] .  \notag
\end{eqnarray}

D'autre part,%
\begin{eqnarray}
&&\left\Vert x+th\right\Vert ^{2^{n-1}}+\left\Vert x-th\right\Vert
^{2^{n-1}}-2\left\Vert x\right\Vert ^{2^{n-1}}  \notag \\
&\leq &\sqrt{2}(\left\Vert x+th\right\Vert ^{2^{n}}+\left\Vert
x-th\right\Vert ^{2^{n}})^{\frac{1}{2}}-2\left\Vert x\right\Vert ^{2^{n-1}}
\label{io} \\
&\leq &\sqrt{2}\left[ (\left\Vert x+th\right\Vert ^{2^{n}}+\left\Vert
x-th\right\Vert ^{2^{n}})^{\frac{1}{2}}-\sqrt{2}\left\Vert x\right\Vert
^{2^{n-1}}\right] .  \notag
\end{eqnarray}

Il est clair que d'apr\`{e}s (\ref{gt}) et (\ref{io}) que lim$_{R\ni
t\rightarrow 0}$ $\frac{\left\Vert x+th\right\Vert ^{2^{n-1}}+\left\Vert
x-th\right\Vert ^{2^{n-1}}-2\left\Vert x\right\Vert ^{2^{n-1}}}{\left\vert
t\right\vert }=0.\blacksquare $

\begin{remark}
\label{C}Dans le lemme pr\'{e}c\'{e}dent on peut remplacer $S_{X}$ par $%
X-\left\{ 0\right\} .$
\end{remark}

\begin{lemma}
\label{UN}Soient ($X,\left\Vert .\right\Vert )$ un espace norm\'{e}, $h,x\in
X$ et $n\in \mathbb{N}$. Alors il existe une constante $C_{n}$ ind\'{e}%
pendante de $h,x$ telle que $\sup_{t\in \left[ -1,1\right] -\left\{
0\right\} }\left\vert \frac{\left\Vert x+th\right\Vert ^{2^{n}}+\left\Vert
x+th\right\Vert ^{2^{n}}-2\left\Vert x\right\Vert ^{2^{n}}}{t}\right\vert
\leq C_{n}(\left\Vert x\right\Vert ^{2^{n}}+\left\Vert h\right\Vert
^{2^{n}}).$
\end{lemma}

D\'{e}monstration: Il suffit de montrer que $\sup_{t\in \left[ -1,1\right]
-\left\{ 0\right\} }\left\vert \frac{\left\Vert x+th\right\Vert
^{2^{n}}-\left\Vert x\right\Vert ^{2^{n}}}{t}\right\vert \leq
C_{n}(\left\Vert x\right\Vert ^{2^{n}}+\left\Vert h\right\Vert ^{2^{n}}),$
car $\sup_{t\in \left[ -1,1\right] -\left\{ 0\right\} }\left\vert \frac{%
\left\Vert x+th\right\Vert ^{2^{n}}+\left\Vert x+th\right\Vert
^{2^{n}}-2\left\Vert x\right\Vert ^{2^{n}}}{t}\right\vert \leq 2\sup_{t\in %
\left[ -1,1\right] -\left\{ 0\right\} }\left\vert \frac{\left\Vert
x+th\right\Vert ^{2^{n}}-\left\Vert x\right\Vert ^{2^{n}}}{t}\right\vert .$
Nous montrons cela par r\'{e}currence sur $n.$ L'\'{e}nonc\'{e} est vrai
pour $n=0.$ Supposons qu'il est vrai pour $n$; montrons qu'il est vrai pour $%
n+1.$ Remarquons que

\begin{eqnarray*}
&&\text{sup}_{t\in \left[ -1,1\right] -\left\{ 0\right\} }\left\vert \frac{%
\bigl\Vert x+th\bigr\Vert^{2^{n+1}}-\bigl\Vert x\bigr\Vert^{2^{n+1}}}{t}%
\right\vert \\
&=&\text{sup}_{t\in \left[ -1,1\right] -\left\{ 0\right\} }\left[ \left\vert 
\frac{\bigl\Vert x+th\bigr\Vert^{2^{n}}-\bigl\Vert x\bigr\Vert^{2^{n}}}{t}%
\right\vert \times \left[ \left\vert \bigl\Vert x+th\bigr\Vert^{2^{n}}+%
\bigl\Vert x\bigr\Vert^{2^{n}}\right\vert \right] \right] \\
&\leq &C_{n}(\bigl\Vert x\bigr\Vert^{2^{n}}+\bigl\Vert h\bigr\Vert%
^{2^{n}})\times (2^{-1}+1)\bigl\Vert x\bigr\Vert^{2^{n}}+2^{-1}\bigl\Vert h%
\bigr\Vert^{2^{n}}).
\end{eqnarray*}

D'autre part, $2(\bigl\Vert x\bigr\Vert^{2^{n}}\times (\bigl\Vert h\bigr\Vert%
^{2^{n}}\leq ((\bigl\Vert x\bigr\Vert^{2^{n+1}}+(\bigl\Vert h\bigr\Vert%
^{2^{n+1}}).$ Il en r\'{e}sulte qu'il existe $C_{n+1}$ telle que

sup$_{t\in \left[ -1,1\right] -\left\{ 0\right\} }\left\vert \frac{%
\bigl\Vert
x+th\bigr\Vert^{2^{n+1}}-\bigl\Vert x\bigr\Vert^{2^{n+1}}}{t}\right\vert
\leq C_{n+1}(\bigl\Vert x\bigr\Vert^{2^{n+1}}+\bigl\Vert h\bigr\Vert%
^{2^{n+1}}).\blacksquare $

\begin{lemma}
\label{UUU}Soient $X$ un espace de Banach et $n\in \mathbb{N}^{\ast };$
suppsons que $X^{\ast }$ est G\^{a}teaux diff\'{e}rentiable. Alors $%
VB^{r}(\Omega ,X^{\ast })$ est G\^{a}teaux diff\'{e}rentiablse, o\`{u} $%
r=2^{n}.$
\end{lemma}

D\'{e}monstration. D'apr\`{e}s le lemme \ref{vvv}, il suffit de montrer que
lim$_{R\ni t\rightarrow 0}\frac{\left\Vert f+th\right\Vert ^{r}+\left\Vert
f-th\right\Vert ^{r}-2\left\Vert f\right\Vert ^{r}}{t}=0,$ pour tout $f\in
VB^{r}(\Omega ,X^{\ast })-\left\{ 0\right\} $ et $h\in VB^{r}(\Omega
,X^{\ast }).$ Soient $f\in VB^{r}(\Omega ,X^{\ast })-\left\{ 0\right\} ,$ $%
t>0$ et $h\in VB^{r}(\Omega ,X^{\ast });$ nous avons que

\begin{eqnarray*}
&&\frac{\left\Vert f+th\right\Vert _{VB^{r}(\Omega ,X^{\ast
})}^{r}+\left\Vert f-th\right\Vert _{VB^{r}(\Omega ,X^{\ast
})}^{r}-2\left\Vert f\right\Vert _{VB^{r}(\Omega ,X^{\ast })}^{r}}{t} \\
&=&\mathop{\displaystyle \int}\limits_{\Omega }\frac{\bigl\Vert f(\omega
)+th(\omega )\bigr\Vert_{X^{\ast }}^{r}+\bigl\Vert f(\omega )-th(\omega )%
\bigr\Vert_{X^{\ast }}^{r}-2\bigl\Vert f(\omega )\bigr\Vert_{X^{\ast }}^{r}}{%
t}d\mu (\omega ) \\
&=&\mathop{\displaystyle \int}\limits_{\left\{ \omega ;\text{ }f(\omega
)=0\right\} }\frac{\bigl\Vert f(\omega )+th(\omega )\bigr\Vert_{X^{\ast
}}^{r}+\bigl\Vert f(\omega )-th(\omega )\bigr\Vert_{X^{\ast }}^{r}-2%
\bigl\Vert f(\omega )\bigr\Vert_{X^{\ast }}^{r}}{t}d\mu (\omega )+ \\
&&\mathop{\displaystyle \int}\limits_{\left\{ \omega ;\text{ }f(\omega )\neq
0\right\} }\frac{\bigl\Vert f(\omega )+th(\omega )\bigr\Vert_{X^{\ast }}^{r}+%
\bigl\Vert f(\omega )-th(\omega )\bigr\Vert_{X^{\ast }}^{r}-2\bigl\Vert %
f(\omega )\bigr\Vert_{X^{\ast }}^{r}}{t}d\mu (\omega ) \\
&=&\mathop{\displaystyle \int}\limits_{\left\{ \omega ;\text{ }f(\omega
)=0\right\} }\frac{2t^{r}\bigl\Vert h(\omega )\bigr\Vert_{X^{\ast }}^{r}}{t}%
d\mu (\omega )+ \\
&&\mathop{\displaystyle \int}\limits_{\left\{ \omega ;\text{ }f(\omega )\neq
0\right\} }\frac{\bigl\Vert f(\omega )+th(\omega )\bigr\Vert_{X^{\ast }}^{r}+%
\bigl\Vert f(\omega )-th(\omega )\bigr\Vert_{X^{\ast }}^{r}-2\bigl\Vert %
f(\omega )\bigr\Vert_{X^{\ast }}^{r}}{t}d\mu (\omega ).
\end{eqnarray*}

Il est evident que lim$_{R\ni t\rightarrow 0}$ $\mathop{\displaystyle \int}%
\limits_{\left\{ \omega ;\text{ }f(\omega )=0\right\} }\frac{2t^{r}%
\bigl\Vert
h(\omega )\bigr\Vert_{X^{\ast }}^{r}}{t}d\mu (\omega )=0.$ D'autre part d'apr%
\`{e}s le lemme \ref{vvv}

lim$_{R\ni t\rightarrow 0}\frac{\bigl\Vert f(\omega )+th(\omega )\bigr\Vert%
_{X^{\ast }}^{r}+\bigl\Vert f(\omega )-th(\omega )\bigr\Vert_{X^{\ast
}}^{r}-2\bigl\Vert f(\omega )\bigr\Vert_{X^{\ast }}^{r}}{t}=0.$ En utilisant
le lemme \ref{UN} et le th\'{e}or\`{e}me de convergence domin\'{e}e, on voit
que%
\begin{equation*}
lim_{R\ni t\rightarrow 0}\mathop{\displaystyle \int}\limits_{\left\{ \omega ;%
\text{ }f(\omega )\neq 0\right\} }\frac{\bigl\Vert f(\omega )+th(\omega )%
\bigr\Vert_{X^{\ast }}^{r}+\bigl\Vert f(\omega )-th(\omega )\bigr\Vert%
_{X^{\ast }}^{r}-2\bigl\Vert f(\omega )\bigr\Vert_{X^{\ast }}^{r}}{t}d\mu
(\omega )=0.
\end{equation*}

Ceci implique que lim$_{R\ni t\rightarrow 0}\mathop{\displaystyle \int}%
\limits_{\Omega }\frac{\bigl\Vert f(\omega )+th(\omega )\bigr\Vert_{X^{\ast
}}^{r}+\bigl\Vert f(\omega )-th(\omega )\bigr\Vert_{X^{\ast }}^{r}-2%
\bigl\Vert f(\omega )\bigr\Vert_{X^{\ast }}^{r}}{t}d\mu (\omega
)=0.\blacksquare $

Dans \cite{Le-Su} on montre que si le Banach X est Fr\'{e}chet diff\'{e}%
rentiable, alors $L^{p}(\Omega ,X)$ a la m\^{e}me propri\'{e}t\'{e}. Dans le
th\'{e}or\`{e}me suivant on montre un r\'{e}sultat analogue lorsque $X^{\ast
}$ est G\^{a}teaux diff\'{e}rentiable.

\begin{theorem}
\label{ts}Soient $(\Omega ,\mu )$ un espace de probabilit\'{e} et $X$ un
espace de Banach. Supposons que $X^{\ast }$ est G\^{a}teaux diff\'{e}%
rentiable. Alors $L^{p}(\Omega ,X)^{\ast }$ est G\^{a}teaux diff\'{e}%
rentiable, pour tout $1<p<+\infty $.
\end{theorem}

D\'{e}monstration: Soit $Y$ un sous espace ferm\'{e} s\'{e}parable de $%
L^{p}(\Omega ,X);$ il existe un sous espace ferm\'{e} $X_{1}$ s\'{e}parable
de $X$ tel que $Y$ se plonge isom\'{e}triquement dans $L^{p}(\Omega ,X_{1}).$
D'apr\`{e}s \cite[Lemme 4.4]{Da} $(iii)\Rightarrow (i)$ appliqu\'{e} \`{a} $%
Z=L^{p}(\Omega ,X_{1}),$ il suffit de montrer que $Z^{\ast }$ est G\^{a}%
teaux diff\'{e}rentiable. Remarquons d'apr\`{e}s \cite[Lemme 4.4]{Da} ($%
(i)\Rightarrow (ii),$ que $(X_{1})^{\ast }$ est G\^{a}teaux diff\'{e}%
rentiable). Donc on peut supposer que $X$ est s\'{e}parable D'apr\`{e}s \cite%
[p.~349]{Blas}, \cite[Chap~II-13-3, Corollary~1]{Din} $(L^{p}(\Omega
,X))^{\ast }=VB^{p^{\prime }}(\Omega ,X^{\ast }),$ o\`{u} $p^{\prime }$ est
le conjugu\'{e} de $p.$ D'autre part, d'apr\`{e}s le lemme \ref{UUU} $%
VB^{r}(\Omega ,X^{\ast })=(L^{r^{\prime }}(\Omega ,X))^{\ast },VB^{\rho
}(\Omega ,X)^{\ast }=(L^{\rho ^{\prime }}(\Omega ,X))^{\ast }$ sont G\^{a}%
teaux diff\'{e}rentiablses, o\`{u} $r=2^{n},\rho =2^{m}.$ Mais d'apr\`{e}s 
\cite[Th.5.1.2]{Ber-Lof} que $(L^{r^{\prime }}(\Omega ,X),L^{\rho ^{\prime
}}(\Omega ,X))_{\theta }=L^{p}(\Omega ,X),$ o\`{u} $\frac{1}{p}=\frac{%
1-\theta }{r^{\prime }}+\frac{\theta }{\rho ^{\prime }}.$ Finalement d'apr%
\`{e}s \cite[Th.4.5.1]{Ber-Lof} $\ (VB^{r}(\Omega ,X^{\ast }),VB^{\rho
}(\Omega ,X^{\ast }))^{\theta }=VB^{p^{\prime }}(\Omega ,X^{\ast }),$ par
cons\'{e}quent $VB^{p^{\prime }}(\Omega ,X^{\ast })=L^{p}(\Omega ,X))^{\ast
} $ est G\^{a}teaux diff\'{e}rentiable d'apr\`{e}s \cite[Th.4.5]{Da1}.$%
\blacksquare $\textit{\ }

\end{document}